\pdfoutput=1 

\documentclass[11pt]{article}
\usepackage{amsmath,amsfonts,amssymb}
\usepackage{fancyhdr,graphicx,amsmath, amscd}
\usepackage{color}

\def\be{\begin{equation}}
\def\ee{\end{equation}}
\def\bea{\begin{eqnarray}}
\def\eea{\end{eqnarray}}
\def\bt{\begin{theorem}}
\def\et{\end{theorem}}
\def\bl{\begin{lemma}}
\def\el{\end{lemma}}
\def\br{\begin{remark}}
\def\er{\end{remark}}
\def\bc{\begin{corollary}}
\def\ec{\end{corollary}}
\def\bd{\begin{definition}}
\def\ed{\end{definition}}
\def\eprop{\end{proposition}}
\def\bprop{\begin{proposition}}
\def\a{\alpha}
\def\b{\beta}

\def\d{\delta}

\def\k{\kappa}

\def\s{\sigma}
\def\t{\tau}

\def\D{\Delta}

\def\cE{\mathcal{E}}

\def\bbR{\mathbb{R}}

\def\b1{B_{1}^}

\def\pt{\partial}

\def\ba{\begin{array}}
\def\ea{\end{array}}
\def\ben{\begin{enumerate}}
\def\een{\end{enumerate}}
\newtheorem{theorem}{Theorem}
\newtheorem{lemma}{Lemma}
\newtheorem{remark}{Remark}
\newtheorem{proposition}{Proposition}
\newtheorem{corollary}{Corollary}

\newtheorem{definition}{Definition}

\begin{document}
\title{On the Homothety Conjecture
}
\author{Elisabeth M. Werner \thanks{Partially supported by an NSF grant, a FRG-NSF grant and  a BSF grant.}
\and Deping Ye \thanks{Support by NSF-FRG DMS: 0652571.}}
\date{ }
\maketitle
\begin{abstract}
Let $K$ be a convex body in $\bbR^n$ and $\delta>0$. The homothety
conjecture asks:  Does $K_{\delta}=c K$ imply that $K$ is an
ellipsoid? Here $K_{\delta}$ is the (convex) floating body and $c$
is a constant depending on $\delta$ only.  In this paper we prove
that the homothety conjecture holds true in the class of  the
convex bodies $B^n_p$, $1\leq p\leq \infty$, the unit balls of
$l_p^n$; namely, we show that $(B^n_p)_{\delta} = c B^n_p$ if and
only if $p=2$. We also show that the homothety conjecture is true
for a general convex body $K$ if $\delta$ is small enough. This
improvs earlier results by Sch\"utt and Werner \cite{SW1994} and
Stancu \cite{Stancu2009}.
\end{abstract}

\section{Introduction}

 Floating bodies  appear in many contexts and have been widely studied
(see e.g. \cite{BaranyVitale1993, BuchtaReitzner1997,
Leichtweiss1986, MeyerReisner1989,  MeyerReisner1991, MW2,
Schutt1991, Schutt1998, SW1990, Stancu2006, Stancu2009,
StancuWerner2009, Werner1999, WY2008}). The homothety conjecture
is among the   problems related to  floating bodies that was open
for a long time. It asks: \vskip 2mm
 {\it Does $K$ have to be an ellipsoid, if $K$ is homothetic to
$K_\d$ for some $\d>0$?}

\vskip 2mm
 In \cite{SW1994}, Sch\"utt and Werner obtained a
(partial) positive solution to this conjecture. They showed that
if there is a sequence $\d_k\rightarrow 0$ such that  $K_{\d_k}$
is homothetic to $K$ for all $k \in \mathbb{N}$ (with respect to
the same center of homothety), then $K$ is an ellipsoid.

\vskip 2mm Stancu \cite{Stancu2009} (see also \cite{Stancu2006})
proved that if for a  convex body $K$ with boundary of class $C^{2
}_+$ there exists a positive number $\d(K)$ such that $K_{\d}$ is
homothetic to $K$ for some $\d<\d(K)$,  then  and only then  $K$
is an ellipsoid.

\vskip 2mm These results  are not  completely satisfactory for
different reasons. The first one requires a sequence $K_{\d_k}$ of
convex bodies to be homothetic to the body $K$. The second one
needs smoothness assumptions  on $\partial K$ and, in addition,
$\d$ has to be sufficiently small, but no estimates are given how
small. Moreover, these results do not work even in the case of
very basic convex bodies, such as $B^n_p$, $1\leq p\leq \infty$,
the unit balls of $l^n_p$.

\vskip 2mm It should also be noted in this context, that Milman
and Pajor \cite{Milman-Pajor1989} showed that if $K$ is a
symmetric convex body, then for large $\d$, the floating body
$K_{\delta }$ is always uniformly, up to a factor $C(\d)$,
isomorphic to an ellipsoid, namely to the dual of the Binet
ellipsoid from classical mechanics. \vskip 3mm In this paper, we
give a positive solution of the homothety conjecture. Namely we
prove \vskip 2mm \noindent {\bf Theorem \ref{Homothety:C2}} {\em
Let $K$ be a convex body in $\mathbb{R}^n$. There exists a
positive number $\d(K)$, such that $K_{\d}$ is homothetic to $K$
for some $\d \leq \d(K)$, if and only if $K$ is an ellipsoid. }

\vskip 2mm \noindent  Our proof is different from Stancu's proof.
No smoothness assumptions are required. In fact, the main
ingredient in our proof is to show that the homothety assumption
implies that $\partial K$ is $C^2_+$. This is done in Lemmas
\ref{C1} and \ref{posdef}.

\vskip 2mm Our proof  of Theorem \ref{Homothety:C2} gives a
possibility to estimate the threshold $\d(K)$ for  convex bodies
$K$ in $\mathbb{R}^n$ that have sufficiently smooth boundary. This
is done in Section 5.

\vskip 2mm  We also show the following theorem, which gives a
positive solution of the homothety conjecture
 for the $B^n_p$ balls, $1\leq p\leq \infty$, and
their affine images without any requirements on the size of
$\delta$.

\vskip 2mm \noindent {\bf Theorem \ref{Homothetic:Unit:ball:lp}}
{\em Let $B^n_p, 1\leq p\leq \infty$ be the unit ball of $l^n_p$.
Let $0<\d<\frac{|B^n_p|}{2}$. Then $(B^n_p)_{\d}=c B^n_p$ for some
$0<c<1$ if and only if $p=2$.}

\vskip 4mm Throughout the paper we will use the following
notations. For $1\leq p \leq \infty$, $B^n_p=\{x  \in\bbR^n:
\|x\|_p  \leq 1\}$ is  the unit ball of $l^n_p= (\bbR^n,  \|  \
\|_p)$ where for $1 \leq p < \infty$
$$
\|x\|_p = \left(\sum _{i=1}^n|x_i|^p \right)^\frac{1}{p},
$$
and for $p=\infty$
$$
\|x\|_\infty = \mbox{max} _{1 \leq i \leq n}|x_i|.
$$
Let $u\in S^{n-1}=\partial B^n_2$, the boundary of the Euclidean
unit ball $B^n_2$. Then $H(x,u)=\{y\in \bbR^n, \langle y, u
\rangle= \langle x,u \rangle \}$ is the hyperplane  through $x$
with outer unit normal vector $u$. The two half-spaces generated
by $H(x,u)$ are $H^-(x,u)=\{y\in \bbR^n, \langle y, u \rangle\geq
\langle x,u\rangle\}$ and $H^+(x,u)=\{y\in \bbR^n, \langle y, u
\rangle\leq \langle x,u\rangle\}$.  $\| \ \|$ is the Euclidean
norm on $\mathbb{R}^n$.

\vskip 2mm For a convex body $K$ in $\mathbb{R}^n$ and $x \in
\partial K$, the boundary of $K$, $N_K(x)$ denotes the outer unit
normal vector to $K$   and $\kappa_K(x)$ the Gauss curvature of
$K$ at $x$. $N_K(x)$ exists almost everywhere
(see\cite{Rock1970}). $\mbox{int}(K)$ is the interior of $K$. We
write $|K|$ or $\mbox{vol}_n(K)$ for the volume of  $K$. We say
that $K$ is  in $C^2_+$, if $\partial K$ is $C^2$ and has
everywhere strictly positive Gauss curvature. Without loss of
generality, we will assume throughout the paper that, unless
specified otherwise, $0$ is the center of gravity of $K$  and that
$K$ homothetic to $K_\delta$ is meant with $0$ as the center of
homothety.

\vskip 4mm
The paper is organized as follows. In Section 2  we provide background and prove
some  properties of the convex floating body that are needed in the next sections.
In Section 3 we prove that the homothety conjecture
holds true for $B_p^n$, $1\leq p\leq \infty$. In Section 4 we
prove Theorem \ref{Homothety:C2}. Moreover, we give (partial)
positive solutions of a generalized homothety conjecture. In Section 5 we give
estimates on the threshold $\d(K)$.


\section{The convex floating body and homothety}
Let $K$ be a convex body with $0\in \mbox{int}(K)$ and $\d$ be a positive
number such that $\d <\frac{|K|}{2}$. The convex floating body is
defined as follows \cite{SW1990}.
\begin{definition}\cite{SW1990} \label{convex:floating:body}
Let $K$ be a convex body in $\bbR^n$.  The convex floating body
$K_{\d}$ is  the intersection of all halfspaces $H^+$ whose
defining hyperplanes $H$ cut off a set of volume at most $\delta$
from $K$,
$$K_{\d}=\bigcap_{|H^-\cap K|\leq \d} {H^+}. $$
\end{definition}
Clearly, $K_0=K$ and $K_{\d}\subset K$ for all $\d \geq 0$.
Moreover, the convex floating body has the following property: for
all (invertible) affine maps $T$ on $\bbR^n$ and for all $\d>0$
\begin{equation}\label{Affine:map:floating:body}
(TK)_{\d}=T\left(K_{\frac{\d} {|det(T)|}} \right).
\end{equation}
Here $|\mbox{det}(T)|$ is the absolute value of the determinant of $T$.
In particular, for an affine map $T$ with $|\mbox{det}(T)|=1$,
$(TK)_{\d}=T(K_{\d})$ for all  $\d\geq0$.

\vskip 3mm An ellipsoid  $\cE$ is the affine image of $B^n_2$,
$\cE=T(B^n_2)$,  for some invertible affine map $T$ on $\bbR^n$.
It is easy to see that $(B^n_2)_{\d}=cB^{n}_2$ for all $\d\geq 0$
and for a constant $c=c(\d)< 1$ depending on $\d$ only. Hence one
gets with (\ref{Affine:map:floating:body}) that for all ellipsoids
$\cE$
$$
\cE_{\d}=(T(B^n_2))_{\d}=T\left(  (B^{n}_2)_{\frac{\d} {|\mbox{det}
(T)| }} \right) = T \left( c\left(\frac{\d} {| \mbox{det} (T)
|}\right)  B^{n}_2 \right)= c\left(\frac{\d} {| \mbox{det} (T)
|}\right) \cE
$$  for all $\d\geq 0$ and for
some constant $c\left(\frac{\d} {| \mbox{det} (T) | }\right)<1$. In
other words, if the homethety conjecture holds true, then
$K_{\d}=cK$ for some $\d>0$ and some constant $0<c<1$ if and only
if $K$ is an ellipsoid in $\bbR^n$.

\vskip 3mm
Now we make some general observations concerning
homothety of $K$ with one of its convex floating bodies. First,
only a strictly convex body can be homothetic to one of its convex
floating bodies. This is a consequence of the flowing lemma proved in
\cite{SW1994}.
\par
\bl \label{strict:SW1994} \cite{SW1994}
Let $K$ be a convex body in
$\mathbb{R}^n$, and let $0<\d < |K|/2$.  Then $K_\d$ is strictly
convex. \el

\par Thus, in particular, no polytope can be homothetic to
one of its floating bodies. The next lemma and its proof is
almost identical to Lemma 3 of \cite{SW1994}. We give the proof
for completeness.
\par
\bl \label{C1}
 Let $K$ be a convex body in $\mathbb{R}^n$, and let $\d >0$.
 If  $K_\d$ is homothetic to $K$, then $\partial K$ is of class $C^1$.
\el

\par
\noindent {\bf Proof.} Suppose that $\partial K$ is not of class
$C^1$. Then there is $x_0 \in \partial K$ so that $\partial K$ has
two different supporting hyperplanes,  $H_1$ and $H_2$,  passing
through $x_0$.
 We may assume
that there are two sequences $(x_k^1)_{k \in \mathbb{N}}$ and
$(x_k^2)_{k \in \mathbb{N}} $ on $\partial K$ converging to $x_0$
so that we have: \\
(i) the supporting hyperplanes $H_k^i$  through $x_k^i$, $k \in \mathbb{ N}$, $i = 1, 2$ are unique; \\
(ii) the hyperplanes  $H_k^i$ converge to $H_i$, $ i = 1, 2$; \\
(iii) $\mbox{lim}_{k\rightarrow \infty} N_K(x_k^i)$  is orthogonal to $H_i$, $ i = 1, 2$.\\

See the proof of Lemma 3 in \cite{SW1994} for the construction of
these sequences.
 We
choose a coordinate system such that $x _0 = 0$, $$H_1=\{x\in
\mathbb{R}^n: x(n)=ax(n - 1)\},$$ and $$H_2=\{x\in \mathbb{R}^n:
x(n)=bx(n - 1)\},$$ with $b < a$ where $x(l)$ denotes the $l$-th
coordinate of the vector $x \in \mathbb{R}^n$. Let $x_\d$ be the
point on $\partial K_\d$ that corresponds to $x_0$ by homothety.
We can assume that in the chosen coordinate system $x_\d=(0,
\dots, 0, x_\d(n))$ with $x_\d(n) >0$. Let
$$
H= \left\{x + x_\d \in \mathbb{R}^n: x(n)=\frac{a+b}{2} x(n -
1)\right\}.
$$ By homothety, one sees that $H$ is a support hyperplane of
$\partial K_{\d}$. For $\alpha, \beta \in \mathbb{R}$, let
$$
M (\alpha, \beta)=\left \{x\in \mathbb{R}^n: x \in K \cap H^-, x(n
- 1) = \alpha, x(n) = \beta \right \},
$$
and let $M^*$ be such a set for which the $(n- 2)$-dimensional
volume is maximal
$$
\mbox{vol}_{n- 2}(M^*) = \mbox{max} _{(\alpha, \beta) \in
\mathbb{R}^2}\mbox{vol}_{n- 2}  \left (M (\alpha, \beta)\right).
$$
We consider the set in the $(x(n- 1), x(n))$-plane that consists
of all points $(\alpha, \beta)$ so that $ M (\alpha, \beta) \neq
\emptyset$. This set is contained in the triangle $T $ bounded by
the lines $x(n)= ax(n- 1)$, $x(n)= bx(n- 1)$, and $x(n)= \frac{a +
b}{2} x(n - 1) + x_\d(n)$. Therefore
\begin{eqnarray}\label{C1.1}
\d \leq |K \cap  H^-| &=& \int_{\mathbb{R}^ 2} M (\alpha, \beta) d
(\alpha, \beta) \leq
 \mbox{vol}_{n- 2}(M^*) \mbox{vol}_{ 2}(T)   \nonumber \\
 & = & \mbox{vol}_{n- 2}(M^*)  \frac{2 |x_\d(n)|^2}{a-b}.
 \end{eqnarray}
Now we consider the sets $K \cap  (x_\d+H_i)^-$,   $ i = 1, 2$. It
follows from (i), (ii) and Lemma 2  of \cite{SW1994} that
\begin{equation*}
 |K \cap  (x_\d+H_i)^- |=  \d,   \ \ \ \\    i = 1, 2.
\end{equation*}
This is because $x_{\d}+H_i, i=1,2$ are tangent hyperplanes of
$K_{\d}$. Let $\varepsilon > 0$ be given. By (iii) we may choose
$x_1$ and $x_2$ in $\partial K$ such that for $i=1,2$
\begin{eqnarray} \label{C1.3}
\bigg(\big|x_i(n-1) - P_{H_i}(x_i)(n -1)\big|^{2} +\big |x_i(n) -
P_{H_i}(x_i)(n )\big|^{2} \bigg)^ \frac{1}{2}
\leq  \nonumber \\
\varepsilon \  \bigg(\big|P_{H_i}(x_i)(n -1)\big|^{2} +
\big|P_{H_i}(x_i)(n )\big|^{2}\bigg) ^ \frac{1}{2},
\end{eqnarray}
where $P_{H_i}$ denotes the orthogonal projection onto $H_ i$,  $i
= 1, 2$. As $K_\delta \subset K$, we can choose $x_1$ and $x_2$ in $\partial K$ such that for $i
=1,2$, $P_{H_i + x_\d}(x_i) \in K$ and (\ref{C1.3}) still holds.
\vskip 2mm
Moreover, as $K$ is strictly convex, by Lemma
\ref{strict:SW1994} and homothety, we have that for $i=1, 2$
\begin{eqnarray*}
\bigg(\big|x_i(n-1) - P_{H_i}(x_i)(n -1)\big|^{2} +\big |x_i(n) -
P_{H_i}(x_i)(n )\big|^{2} \bigg)^ \frac{1}{2}  > 0.
\end{eqnarray*}
It follows that there is a constant $c>0$ that depends only on $K$
so that
\begin{equation}\label{xdelta}
x_\d(n) \leq c \bigg(\big|x_i(n-1) - P_{H_i}(x_i)(n -1)\big|^{2}
+\big |x_i(n) - P_{H_i}(x_i)(n )\big|^{2} \bigg)^ \frac{1}{2}.
\end{equation}
For $i=1$ or $i=2$ the convex hull
$$
\big[M^*, [x_i, P_{H_i + x_\d}(x_i)]\big]\subset K \cap
(x_\d+H_i)^-.
$$
We may assume that $i=1$. Then
\begin{eqnarray} \label{C1.4}
\d=  |K \cap  (x_\d+H_1)^- | &\geq & \bigg|  \big[M^*, [x_1, P_{H_1 + x_\d}(x_i)]\big] \bigg| \nonumber \\
&\geq  & \frac{1}{n^2} \mbox{vo}l_{n-2} (M^*) \ d(x_1, H_{M^*}) \
d\left(P_{H_1 + x_\d}(x_1), \tilde{H}\right),
\end{eqnarray}
where $H_{M^*}$ is the $(n-2)$-dimensional flat containing $M^*$
and $\tilde{H}$ is the plane containing $M^*$ and $x_1$. If
$\varepsilon$  is sufficiently small, we obtain by elementary
computation
$$
 d(x_1, H_{M^*}) \geq c_1 \left(\big|P_{H_1}(x_1)(n -1)\big|^{2} + \big|P_{H_1}(x_1)(n )\big|^{2} \right)^\frac{1}{2},
$$
and
$$
d\left(P_{H_1 + x_\d}(x_1), \tilde{H}\right) \geq c_2 |x_\d(n)|,
$$
where $c_1, c_2$  depend only on $K$. Therefore we obtain by
(\ref{C1.1}) and (\ref{C1.4})
$$
\frac{2 |x_\d(n)| }{a-b} \geq \frac{c_1\ c_2}{n^2}
\left(\big|P_{H_1}(x_1)(n -1)\big|^{2} + \big|P_{H_1}(x_1)(n
)\big|^{2} \right)^\frac{1}{2}.
$$
By (\ref{C1.3}) we get with a new constant $c > 0$
$$
|x_\d(n)| \geq \frac{c}{\varepsilon} \  \bigg(\big|x_1(n-1) -
P_{H_1}(x_1)(n -1)\big|^{2} +\big |x_1(n) - P_{H_1}(x_1)(n
)\big|^{2} \bigg)^ \frac{1}{2}.
$$
If we choose $\varepsilon$ sufficiently small we get a
contradiction to (\ref{xdelta}).

\vskip 3mm \noindent
\begin{remark}  \label{remark1}
If $K$ is strictly convex and $C^1$, $K_\delta$ need not be $C^1$.
\end{remark}
To see that, consider a ``rounded simplex"  $S$ in $\mathbb{R}^2$
whose vertices have been ``rounded" by putting a small
$\varepsilon B^n_2$ ball at each vertex and whose edges have been
made strictly convex by replacing them by arcs of Euclidean balls
$R B^n_2$ with $R$ very large. Then there are   points of
non-differentiability on that parts of $\partial S_\delta$ that
face  the ``edges" of $S$. Hence homothety is crucial in the
previous lemma.
\par
Note however that if $K$ is strictly convex, $C^1$ and,  in
addition,  symmetric, then $K_\delta$ is $C^2_+$. This was shown
in \cite{MeyerReisner1991}.

\vskip 4mm
The following symmetric matrix is closely related to
the curvature of the floating body \cite{Leichtweiss1986} (see
also \cite{SW1994}). For $x\in \mbox{int}(K)$, $\xi \in S^{n-1}$,
and an orthonormal coordinate system in the plane $H(x, \xi)$ with
origin $x$, let $l(\eta)$ be the line through $x$ with direction
$\eta $ and $y=l(\eta)\cap
\partial K$. Let $r(\eta)=\|x-y\|$ and $\beta (\eta)$ be the
angle between $l(\eta)$ and a tangent line to $\partial K$ through
$y$ whose orthogonal projection onto $H(x, \xi)$ is $l(\eta)$.
Define for $1 \leq i, j \leq n-1$
\begin{equation} \label{Curvature:Matrix:Q}
Q(i,j)=\frac{1}{|K\cap H(x, \xi)|} \int _{S^{n-2}} \eta_i \eta_j \
r^n(\eta)\  \mbox{cot} (\beta(\eta)) \  \,d\s_{n-2}(\eta),
\end{equation}
 where $\s _{n-2}$ is the surface measure on $S^{n-2}$. For $x_\d \in \partial K_\d$, we will
 use $H(x_\d, \d)$ to denote the hyperplane through $x_\d$  cutting off a set of volume $\d$ from
$K$. Such a hyperplane always exists by Lemma 2 of \cite{SW1994}.
\vskip 3mm
Then we have the following lemma
\cite{Leichtweiss1986} (see also \cite {SW1994})
\bl \label{C2}  \cite{Leichtweiss1986} Let $K$ be a convex body  of class $C^1$ and let
$\d>0$. Suppose that for every $x_\d  \in
\partial K_\d$ and every  $H(x_\d, \d )$ the matrix $Q$ is
positive definite. Then $\partial K_{\d} $ is of class $C^2$.
 Moreover, the Gauss curvature of $K_{\d}$ at
$x_{\d} \in
\partial K_{\d}$ can be calculated by
$$
\frac{1}{\k_{K_{\d}}(x_{\d}) } = det(Q),
$$
and hence, $K_\d$ is of class $C^2_+$.
\el

\vskip 3mm
Also the next lemma and its proof is
almost identical to Lemma 6 of \cite{SW1994}. Again, we give the proof
for completeness. We will also use the following:
\par
For every $x \in  \partial K$,  we define $\rho(x)$  to be the radius of the largest Euclidean
ball that is contained in $K$  and whose center lies on the line through $x$
with direction $N_K(x)$,
$$
\{x + t \ N_K(x) : t \in \mathbb{R}\}.
$$
$\rho(x)$ is well defined because $\partial K$  is of class $C^1$.  Let  $ t(x) \geq \rho(x) $  be such
that $x - t (x) \ N_K(x)$  is a center of a ball  with radius $\rho(x)$  that is contained
in $ K$. $\rho(x)$  is a continuous, strictly positive function on $\partial K$  because  $\partial K$ is of
class $C^1$. By compactness there is  $\rho_0 > 0$  so that we have for all $ x \in \partial K$
\begin{equation}\label{rho0}
0 < \rho_0 \leq  \rho(x).
\end{equation}
By Lemma 5 of \cite{SW1994}  and the homothety of $K$  and  $K_\delta$,   we get that there is $R > 0 $ so
that we have for all $x  \in \partial K$
\begin{equation}\label{R}
K \subset  B^n_2\left(x - R \ N_K(x), R\right).
\end{equation}

\bl \label{posdef} If $K$ is of class $C^1$  and
homothetic to $K_\d$ for some $\d \leq \delta_0$, then   the matrix $Q$ is
positive definite for every $x_\d$ in $\partial K_\d$.
Moreover $\delta_0$ can be chosen to be
\begin{equation}\label{delta:0:0}
\delta_0=\frac{ \rho_0^{n-1} R  |B^{n-1}_2| }{n\ 2^{n-1} }
\left(1-\left(1- (\frac{\rho_0}{4R})^2\right)^\frac{1}{2}
\right)^n, \end{equation}  where $\rho_0$ and $R$ are as in
(\ref{rho0}) and (\ref{R}). \el

\par \noindent
{\bf Proof.} We want to show that we have for $\delta$
sufficiently small: For all $x \in \partial K_\delta$  and for all
$\eta$ we  have $\mbox{cot} (\beta(\eta)) > 0$. This implies that
(\ref{Curvature:Matrix:Q}) is a positive-definite matrix and by
Lemma \ref{C2} we get that $\partial K_\delta$ is of class  $C^2$
and therefore $\partial K$  is  also of class $C^2$ by homothety.

\vskip 3mm  Let $\rho_0$ and $R$ be as in (\ref{rho0}) and
(\ref{R}). We choose $s > 0$ so that the orthogonal projection of
$$
B^n_2\left(x - R \ N_K(x), R\right) \cap H\left(x - s\ N_K(x),
N_K(x)\right),
$$
is contained in
\begin{equation*}
B^n_2\left(x - t(x) \ N_K(x), \frac{\rho_0}{4} \right) \cap H\left(x-t(x) \ N_K(x), N_K(x)\right).
\end{equation*}
It is easy to see that this holds if  $s=R \left(1-\left(1- (\frac{\rho_0}{4R})^2\right)^\frac{1}{2} \right)$.

\begin{figure}[htp]
\centering
\includegraphics[totalheight=0.45\textheight]{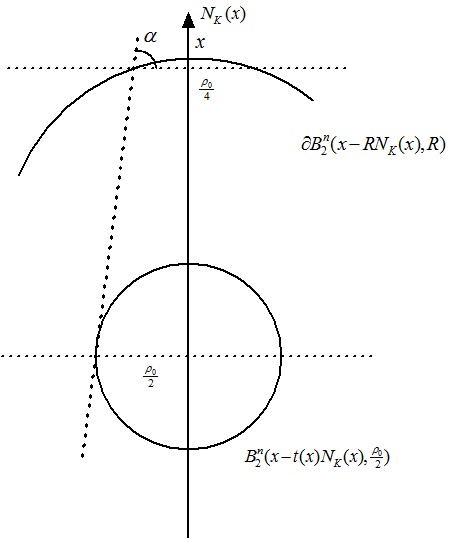}
\caption {}\label{floating:figure:1}
\end{figure}

Now we choose  $\delta_0$. We choose $\delta_0$  so small that  we
have for all $ x \in \partial K$
\begin{equation} \label{cont}
\left|H^-\left(x - s\ N_K(x), N_K(x)\right) \cap  K\right| \geq \delta_0.
\end{equation}

\begin{figure}[htp]
\centering
\includegraphics[totalheight=0.45\textheight]{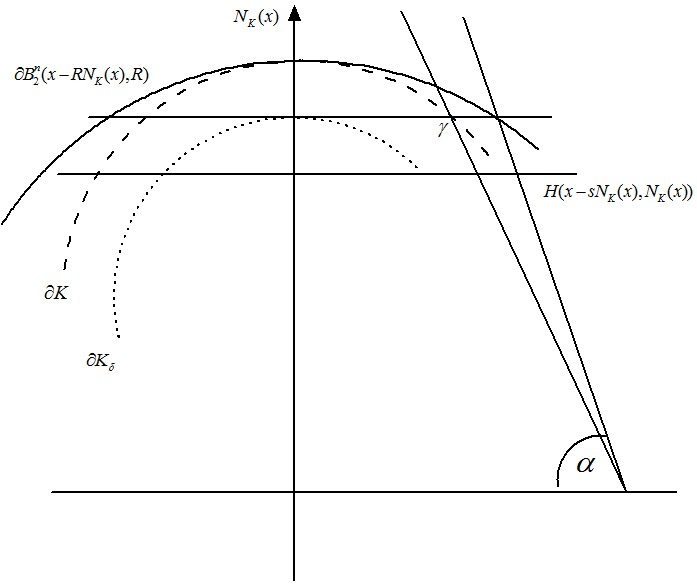}
\caption {}\label{floating:figure:2}
\end{figure}

\noindent From Figure (\ref{floating:figure:1}), it follows that
$$
\mbox{cot} (\alpha)  \geq \frac{\frac{ \rho_0}{2}- \frac{\rho_0}{4}}{t(x)} > 0.
$$
From Figure (\ref{floating:figure:2}), it follows that $
\beta(\eta) \leq \gamma \leq \alpha$ and
$$
\mbox{cot}(\beta(\eta))  \geq \mbox{cot}(\gamma)  \geq  \mbox{cot}(\alpha)  \geq  \frac{\rho_0}{4\  t(x)}.
$$

\noindent We determine now $\delta_0$ from (\ref{cont}) as follows
$$ \delta_0=\frac{ \rho_0^{n-1} R  |B^{n-1}_2| }{n\ 2^{n-1} }
\left(1-\left(1- (\frac{\rho_0}{4R})^2\right)^\frac{1}{2}
\right)^n.
$$
Then for all $\delta \leq \delta_0$, (\ref{cont})  holds.

Indeed, for $x \in \partial K$, let $C\left(t(x), \rho_0\right) = \mbox{conv}\left[ x,  B^{n-1}_2(x-t(x)N_K(x), \rho_0)\right] $ be the cone
with tip at $x$, height $t(x)$ and with the $(n-1)$-dimensional ball centered at $x-t(x)N_K(x)$ and radius $\rho_0$ as base.
Then, as  $B^{n}_2(x-t(x)N_K(x), \rho_0) \subset K$,
\begin{eqnarray*}
 \left|H^-\left(x - s\ N_K(x), N_K(x)\right) \cap  K\right| &\geq&  \left|H^-\left(x - s\ N_K(x), N_K(x)\right) \cap  C\left(t(x), \rho_0\right)\right| \\
&=&
\frac{ \rho_0^{n-1} s^n  |B^{n-1}_2| }{n \  t(x)^{n-1}}   \geq \frac{ \rho_0^{n-1} s^n  |B^{n-1}_2|}{n\ 2^{n-1} R^{n-1}} \\
&=&\frac{ \rho_0^{n-1} R  |B^{n-1}_2| }{n\ 2^{n-1} } \left(1-\left(1- (\frac{\rho_0}{4R})^2\right)^\frac{1}{2} \right)^n \\
&=& \delta_0 \geq \delta,
\end{eqnarray*}
which shows that  (\ref{cont}) holds.

\vskip 3mm \noindent \begin{remark}  \label{remark2}
\par
\noindent
(i) Thus, by Lemmas \ref{strict:SW1994}, \ref{C1}, \ref{C2}, and
\ref{posdef}, a convex body $K$ in $\mathbb{R}^n$ can only be
homothetic to one of its floating bodies $K_\delta$ for small enough $\delta$, if $K$ is strictly convex and of class $C^2_+$.
\par
\noindent
(ii) The example of Remark 1 shows that the assumption of homothety is also crucial
in Lemma \ref{posdef}.
\end{remark}


\section{Homethety conjecture for $K=B^n_p$}
In this section, we show that the homothety conjecture holds true
in the class of  the bodies $B^n_p$ for all $1\leq p\leq
\infty$ and their affine images.
\bt \label{Homothetic:Unit:ball:lp}
Let $B^n_p, 1\leq p\leq \infty$
be the unit ball of $l^n_p$. Let $0<\d<\frac{|B^n_p|}{2}$. Then
$(B^n_p)_{\d}=c\  B^n_p$ for some $0<c<1$ if and only if  $p=2$.
\et

\vskip 3mm  \noindent {\bf Remark. } By
(\ref{Affine:map:floating:body}), the same holds true for  affine
images $T\left(B^n_p\right)$ of $ B^n_p$ under an invertible linar
map $T$ on $\bbR^n$: {\em Let $K = T\left( B^n_p\right)$,  $1\leq
p\leq \infty$. Let $0<\d<\frac{|K|}{2}$ be a constant. $K_{\d}=c\
K$ for some constant $0<c<1$ if and only if $K$ is an ellipsoid.}

\vskip 3mm \noindent {\bf Proof of Theorem
\ref{Homothetic:Unit:ball:lp}.} It was noted above that
$(B^n_2)_{\d}$ is homothetic to $B^n_2$ for all $\d$ with
$0<\d<\frac{|B^n_2|}{2}$.

\vskip 2mm  For the reverse implication, it is enough to consider
$p$ with $1< p <\infty$. Indeed, $B^n_1$ and $B^n_\infty$ are
polytopes and it was observed  above that polytopes cannot be
homothetic to any of their floating bodies.

 \vskip 2mm  We first consider the case when  $1< p < 2$. Then
$\partial B^n_p$ is of class $C^1$ but not of class $C^2$: at $e_n
=(0, \cdots, 0, 1)$, $\partial B^n_p$ is not $C^2$. If $B^n_p$
were homothetic to $(B^n_p)_{\d}$  for some $\delta$, then Lemma
\ref{C2} would imply  that  $ (B^n_p)_\delta$ is $C^2$. Indeed,
for all $x \in \partial (B^n_p)_\delta$, $0 < \beta(\eta) <
\frac{\pi}{2}$ for all $\eta \in S^{n-2}$ and thus $\mbox{cot}
(\beta(\eta)) >0$ or all $\eta \in S^{n-2}$. Therefore the matrix
$Q$ of (\ref{Curvature:Matrix:Q}) is positive definite and  by
Lemma \ref{C2} $(B^n_p)_\d$ and thus, by homothety, $B^n_p$ is
$C^2$, a contradiction.

\vskip 4mm Now we  consider the case $2 \leq p<\infty$. Then
$\partial B^n_p$ is of class $C^2$ and thus, assuming that $B^n_p$
is homothetic to $\left(B^n_p \right)_\d $ for some $\d >0$,
$\left(B^n_p \right)_\d $ is $C^2$.

\vskip 2mm It was shown in \cite{SW5}, that the curvature
$\kappa_{B^n_p}$ at $x \in \partial B^n_p$ is
\begin{equation}
\kappa_{B^n_p}(x) = \frac{(p-1)^{n-1} \prod_{i=1}^n
x_i^{p-2}}{\left(\sum_{i=1}^n |x_i|^{2(p-1)}\right)^\frac{1}{2}}.
\end{equation}
Thus for  $e_n =(0, \cdots, 0, 1)$, $\kappa_{B^n_p}(e_n) = 0$, if
$p >2$ and $\kappa_p(e_n) =1$, if $p=2$. Consequently, as we
assume that $ \left(B^n_p \right)_\d = c B^n_p$ for some $\d >0$,
the curvature at $c  e_n \in \partial \left(B^n_p \right)_\d$ is
different from $0$ only when $p=2$ and for all $p>2$
\begin{equation}\label{zero}
\k_{\left(B^n_p \right)_\d}(c  e_n) =0.
\end{equation}
By Lemma \ref{C2},
$$
\k_{\left(B^n_p \right)_\d}(c  e_n)  = \frac{1}{\mbox{det}(Q)},
$$
which is $0$ if and only if $\mbox{det}(Q) = \infty$.

\vskip 2mm By (\ref{Curvature:Matrix:Q}), the matrix $Q$ has
entries
\begin{equation*}
Q(i,j) = \frac{1}{|K\cap H(c e_n, e_n)|} \int _{S^{n-1}} \eta_i
\eta_j \ r^n(\eta)\  \mbox{cot} (\beta(\eta)) \
\,d\s_{n-2}(\eta),
\end{equation*}
$1 \leq i, j \leq n-1$. Again,  for all $\eta \in S^{n-2}$,  $0 <
\beta(\eta) < \frac{\pi}{2}$. Therefore, $\mbox{det} (Q) <\infty$
and thus the curvature at $c  e_n \in \partial \left(B^n_p
\right)_\d$ is strictly positive which  contradicts (\ref{zero}).

\vskip 2mm We point out that $(B^n_p)_{\d} \in C^2_+$ for all
$p\in (1, \infty)$ and all $\d >0$ can be obtained from results in
\cite{MeyerReisner1991}.

 \vskip 3mm \noindent More generally, in the same way
as Theorem \ref{Homothetic:Unit:ball:lp}, one can prove \bprop
\label{general}
 Let $K$  be a convex body in $\mathbb{R}^n$. If $K$  has a point on the boundary where the Gauss curvature is
either 0 or $\infty$, then K is not homothetic to $K_\delta$ for
any $\d \leq \delta_0$ with $\delta _0$ given as in Lemma
\ref{posdef}.  \eprop

 \noindent {\bf Proof.}
 Let $x_0$ in $\partial K$ be such that
$\kappa_K(x_0)=\infty$ and suppose that $K$ is homothetic to
$K_\d$ for some $\d\leq \delta_0$. Then, by Remark \ref{remark2}, $K$ is
strictly convex and in $C^2_+$. By Lemma \ref{C2},
$$
\k_{K_\d}(c  x_0)  = \frac{1}{\mbox{det}(Q)},
$$
where $c x_0$ is the point on $\partial K_\d$ corresponding to
$x_0$ by homothety. This  is $\infty$ if and only if
$\mbox{det}(Q) = 0$. But by Lemmas \ref{C2} and   \ref{posdef} (
respectively Remark \ref{remark2}),  $K$ is in $C^2_+$, thus $Q$
is positive definite and thus $\mbox{det}(Q) >0$, a contradiction.

\vskip 2mm The case  $\kappa_K(x_0)=0$ is  treated similarly.


\section{Homothety conjecture for general $K$}
In Section 2  we proved the homothety conjecture for $B^n_p, 1\leq
p\leq \infty $ and their affine images. The proof uses the fact
that one only needs to examine one properly chosen direction in
order to be able  to conclude. In this section, we will use two
directions  to prove positive solutions to the homothety
conjecture. Moreover, our approach is robust and can be used to
obtain (partial) positive solutions for {\it generalized homothety
conjecture}.

\bt \label{Homothety:C2}  Let $K $ be a convex body in
$\mathbb{R}^n$. There exists a positive number $\d(K)$, such that
$K_{\d}$ is homothetic to $K$ for some $\d<\d(K)$, if and only if
$K$ is an ellipsoid.
\et

In the next section we provide estimates for the threshold $\d(K)$.

\vskip 2mm  To prove this theorem, we need the following lemmas.
Lemma \ref{Petty:Ellipsoid} was proved in \cite{Petty1985} and
Lemma \ref{Floating:SW1990} was proved in \cite{SW1990}. See also
\cite{WY2008} for similar results.

\bl \label{Petty:Ellipsoid} \cite{Petty1985} Let $K\in C^2_+$ be a convex body with
boundary of class $C^{2}$ and everywhere strictly positive
Gaussian curvature. Let $c(K,n)$ be a constant only depending on
$K$ and $n$. Then $$\k_K(x) \langle x, N_K(x) \rangle ^{-(n+1)}=
c(K,n),\ \ \ \forall x\in \partial K$$ holds true if and only if
$K$ is an ellipsoid.  \el

\bl \label{Floating:SW1990}  \cite{SW1990} Let $K\in C^2_+$ be a convex body with
boundary of class $C^{2}$ and everywhere strictly positive
Gaussian curvature. Then, for any $x\in \pt K$,
\begin{eqnarray*}
\lim _{\d \rightarrow 0} c_n \frac{\langle x, N_{K}(x)\rangle }{n\
\d^{\frac{2}{n+1}}}\left[1-\left(\frac{\|x_{\d}\|}{\|x\|}\right)^n\right]
= \big(\k _{K} (x)\big)^{\frac{1}{n+1}},
 \end{eqnarray*}
 where  $x_{\d}\in \pt (K_{\d}) \cap [0,x]$ and $c_n=2\left(\frac{|B_2^{n-1}|}{n+1}\right)^{\frac{2}{n+1}}$.
 \el

\vskip 3mm \noindent {\bf Proof of Theorem \ref{Homothety:C2}.}
Suppose that $K$ is homothetic to $K_\d$ for some $\d < \d(K)$ with $\delta(K) \leq \delta_0$, where $\delta_0$ is
given by Lemma \ref{posdef}
($\d(K)$ will be determined more precisely later).  Suppose that $K$ not an
ellipsoid. By Lemma \ref{strict:SW1994} and homothety, $K$ is
strictly convex. By Lemma \ref{C1} and homothety, $K$ is $C^1$. By
Lemma \ref{posdef}, homothety and Lemma \ref{C2}, $K$ is in
$C^2_+$. $K\in C^2_+$ implies that $\langle x, N_K(x) \rangle $
and $\k_K(x)$ both are continuous functions on $\partial K$. We
define points $x_M\in \partial K$ and $x_m \in \partial K$  by
\begin{eqnarray*}
(T_M)^{n+1}&=&\frac{\k _K(x_M)}{\langle x_M, N_K(x_M) \rangle
^{n+1}}=\max _{x\in \partial K} \ \frac{\k _K(x)}{\langle x,
N_K(x) \rangle ^{n+1}}, \\
(T_m) ^{n+1} &=&\frac{\k _K(x_m)}{\langle x_m, N_K(x_m) \rangle
^{n+1}} =\min _{x\in \partial K}\  \frac{\k _K(x)}{\langle x,
N_K(x) \rangle ^{n+1}}.
\end{eqnarray*}
As $K$ is not an ellipsoid, Lemma \ref{Petty:Ellipsoid} implies
that $$\frac{\k _K(x_M)}{\langle x_M, N_K(x_M) \rangle
^{n+1}}>\frac{\k _K(x_m)}{\langle x_m, N_K(x_m) \rangle ^{n+1}}.
$$ Equivalently, one has $$\t = \frac{T_M}{T_m} >1. $$

\par
\noindent
By Lemma \ref{Floating:SW1990}, for
$\varepsilon_1=\frac{\tau-1}{3\tau}$, there exists  $\d_1(K)
> 0$, such that, for all $\d \leq \d_1(K)$,
\begin{eqnarray}\label{xM:estimate}
\frac{c_n }{n\ \d^{\frac{2}{n+1}}}\left[1-\left(\frac{\|x_{M,
\d}\|}{\|x_M\|}\right)^n\right] & \geq & \frac{\big(\k _{K}
(x_M)\big)^{\frac{1}{n+1}}}{\langle x_M,
N_{K}(x_M)\rangle} \   (1-\varepsilon_1) \nonumber  \\
&= &T_M \  (1-\varepsilon_1) =T_M \ \left(\frac{2\t
+1}{3\t}\right),
\end{eqnarray}
where $x_{M, \d}=\partial K _{\d} \cap [0, x_M]$.

\vskip 3mm \par Again by Lemma \ref{Floating:SW1990}, for
$\varepsilon_2=\frac{\tau-1}{2}$, there exists $\d_2(K)>0$, such
that, for all $\d \leq \d_2(K)$,
\begin{eqnarray}\label{xm:estimate}
\frac{c_n }{n\ \d^{\frac{2}{n+1}}}\left[1-\left(\frac{\|x_{m,
\d}\|}{\|x_m\|}\right)^n\right] &\leq&  \frac{\big(\k _{K}
(x_m)\big)^{\frac{1}{n+1}}}{\langle x_m, N_{K}(x_m)\rangle} \ (1+
\varepsilon_2) \nonumber \\ & =& T_m (1+ \varepsilon_2) = T_m \
\left(\frac{\t +1}{2}\right),
\end{eqnarray}
where $x_{m,\d}=\partial K _{\d} \cap [0, x_m]$.

\vskip 3mm \par Let $\d (K)=\min\{ \delta_0, \d_1(K), \d_2(K)\}$. If
$K_{\d}=c K$ for some $0<\d<\d(K)$ and for some constant
$c(\d)>0$, then formulas (\ref{xM:estimate}) and
(\ref{xm:estimate}) imply that $$T_M \ \left(\frac{2\t
+1}{3\t}\right) \leq  \frac{c_n \ (1-c^n)}{n\ \d^{\frac{2}{n+1}}}
\leq T_m \ \left(\frac{\t +1}{2}\right),
$$
and hence
$$
T_M \ \left(\frac{2\t +1}{3\t}\right) \leq T_m \ \left(\frac{\t
+1}{2}\right).
$$ Equivalently, $\t\leq 1$ which is a contradiction. Hence $K$ must be an ellipsoid.
\vskip 3mm
\noindent
{\bf Remark.}
We can replace $x_M$ and $x_m$
by any two points $x, y\in \partial K$ such that
$$ \frac{\k _K(x)}{\langle x, N_K(x) \rangle^{n+1}} \neq \frac{\k _K(y)}{\langle y, N_K(y)
\rangle^{n+1}}.
$$ Following the above proof,  one gets
\bt \label{Homothety:C2:1}
Let $K$ be a convex body in $\bbR^n$.
Suppose that there are two points $x, y\in \partial K$, such that,
both $\kappa _K(x)$ and $\kappa _K(y)$ exist and are finite,  and
$$
\k_K(x) \langle x, N_K(x) \rangle ^{-(n+1)}\neq \k_K(y) \langle y, N_K(y) \rangle
^{-(n+1)}.
$$
Then there is a constant $\d(K, x,y)$  depending on $K,
x$ and $y$ such that $K_{\d}$ is not homothetic to $K$ for all $\d
\leq \d (K,x,y)$.
\et

\vskip 4mm
 Analogously, we can ask the following {\it
generalized homothety conjecture}.

\vskip 3mm \noindent {\bf Generalized Homothety Conjecture:} {\it
Let $K$ be a convex body in $\bbR^n$. Does $K_{s}=c K$ for some
$0<s$ and some $c>0$ imply that $K$ is an ellipsoid? Here
$\{K_{s}\}_{s\geq 0}$ is a family of convex bodies constructed
from $K$. }

\vskip 3mm
 \noindent
Besides the convex floating body $K_{\d},$ examples of such  $K_s$
include
\begin{enumerate}
\item The illumination body $K^{\d}$ \cite{Werner1994}
 $$
 K^{\d}=\{x\in \bbR^n: |\mbox{conv}(x,K)|-|K|\leq \d\}.
 $$
\item The convolution body $C(K,t)$ \cite{Kiener1986, Schm1992}
which, for a symmetric  convex body $K$ in $\bbR^n$ and $t\geq 0$
is defined by
$$
C(K, t)=\left\{\frac{x}{2}\in \bbR^n: |K\cap (K+x)|\geq 2t
\right\}.
$$
\item The Santal\'{o}-regions $S(K,t)$ \cite{MW1998}
$$
S(K,t)=\left\{x\in K: |K^x|\leq \frac{1}{t}\right\},
$$
where  $K^x=(K-x)^\circ =\{z\in \bbR^n: \langle z, y-x\rangle \leq
1, \forall y\in K\}$ is the polar body of $K$ with respect to
$x\in K$.
\end{enumerate}
We refer  to \cite{Werner1999, Werner2002, Werner2007} for more general constructions.

\vskip 3mm
The following theorem provides (partial) positive
solutions of the generalized homothety conjecture. Theorem \ref{gen} (i) was proved with a different method in \cite{Stancu2009}.
\bt \label{gen}
Let $K$ be a convex body in
$C^2_+$.
\par
\noindent
(i) \cite{Stancu2009} There exists a positive number $\tilde{\d} (K)$
such that  $K^{\d}$ is homothetic to $K$ for some
$\d<\tilde{\d}(K)$, if and only if $K$ is an ellipsoid.
 \par
 \noindent
 (ii) There exists a positive number
$t(K)$ such that  $C(K, t)$ is homothetic to $K$ for some
$t<t(K)$, if and only if $K$ is an ellipsoid.
 \par
 \noindent
 (iii) There exists a positive number
$\tilde{t} (K)$  such that  $S(K,t)$ is homothetic to $K$ for some
$t<\tilde{t}(K)$, if and only if $K$ is an ellipsoid.
 \et

\vskip 3mm \noindent  {\bf Remark.}  The proof of this theorem is
same as the proof of Theorem \ref{Homothety:C2}. The proof of  (i)
also relies on Lemma 3 in \cite{Werner1994}. For the proof of
(ii), we refer to  results similar to Lemma \ref{Floating:SW1990}
in \cite{Schm1992}. For  (iii), one uses Lemma 13 in
\cite{MW1998}.

\vskip 2mm  Estimates on the thresholds $\tilde{\delta}(K)$,
$t(K)$ and $\tilde{t} (K)$ can be obtained similar to the one for
$\delta(K)$. This  is treated in the next section.

\section {Estimates on the threshold $\delta(K)$}

Our proof  of Theorem \ref{Homothety:C2} gives a possibility to
estimate the threshold $\d(K)$ for a convex body $K$ in
$\mathbb{R}^n$. Thus we assume that $K_\delta$ is homothetic to
$K$ for some $\delta \leq \delta(K)$ with $0$ as the center of
homothety.  Let $x_M\in \partial K$, $x_m \in \partial K$,  $T_M$,
$T_m$ and $\tau$ be as in the previous section. Note that the
points $x_m, x_M$ may not be uniquely determined. We just choose
any two points satisfying the condition. Let
\begin{equation*}
r_m=\kappa_K(x_m)^{-\frac{1}{n-1}} \hskip 10 mm \mbox{and} \hskip
10mm r_M=\kappa_K(x_M)^{-\frac{1}{n-1}}.
\end{equation*}
Let
\begin{equation}\label{epsilon}
a = \min \left\{ 1-\left(\frac{2}{1+\tau}\right)^\frac{n+1}{n-1},
\left(\frac{3 \tau}{1+2 \tau}\right)^\frac{n+1}{n-1}-1 \right\}.
\end{equation}

\vskip 3mm We show \bt Let $K$ be a convex body in $\mathbb{R}^n$
with everywhere on $\partial K$  strictly positive Gauss curvature
and such that $\partial K$ is $C^{3}$. Let $a$ be as in
(\ref{epsilon}).  Then $\delta(K)$ of Theorem \ref{Homothety:C2}
can be chosen to be
\begin{eqnarray*} \d(K)&=&\min\bigg\{\delta_0,  \delta_1, \delta_2, \d _{m},  \d _{M},  \
\frac{(1-a)^{n}\ r_m^n\ |B^n_2| }{2},  \frac{(1-a)^{n}\ r_M^n \
|B^n_2|}{2}\bigg\},
\end{eqnarray*} where $\d_0$ is as in (\ref{delta:0:0}), and the
expressions of $\delta_1, \delta_2, \d _{m},  \d _{M} $ are in the
proof.  \et \vskip 3mm \noindent {\bf Proof.} For
$\alpha=(\alpha_1, \dots \alpha_m) \in \mathbb{N}^m$,  and  an
$|\alpha |$-times continuously differentiable function $g$, let
\begin{eqnarray*}
&&|\alpha |= \alpha_1 + \dots +\alpha_m, \\
&&\alpha != \alpha_1! \dots \alpha_m !, \\
&&D^\alpha g = D_1^{\alpha_1} \dots  D_n^{\alpha_m} g =
\frac{\partial ^{|\alpha|} g}{\partial t_1^{\alpha_1} \dots
\partial t_n^{\alpha_m}},
\end{eqnarray*}
where $D_i^{\alpha_i}=D_j \dots D_j$ is the $\alpha_i$ times
product of $D_i$.

\vskip 3mm  As determining $\delta(K)$ is invariant under affine
transformations of determinant $1$, we can assume that the
ellipsoid approximating  $\partial K$  at $x_m$ is a Eudlidean
ball  and then have (see \cite{SW4}): For $a>0$ given as in
(\ref{epsilon}) above, there exists $\Delta_{a, m}$ such that for
all $\Delta \leq \Delta_{a, m}$
\begin{eqnarray}\label{approx:1}
B^n_2\left(x_m-\bar{r}_m N_K(x_m), \bar{r}_m \right)  \cap H^-\left(x_m- \Delta N_K(x_m), N_K(x_m) \right)  \nonumber \\
 \subseteq K \cap H^-\left(x_m- \Delta N_K(x_m), N_K(x_m) \right),
\end{eqnarray} where $\bar{r}_m=(1-a)r_m$. In addition, we also choose $\Delta_{a, m} < \bar{r}_m$.

\vskip 2mm Assume now that $x_m=0$, that $N_K(x_m)=-e_n$ and  that
the other $(n-1)$ axes of the approximating ellipsoid coincide
with the remaining $(n-1)$ coordinate axes. Locally  we can then
describe $\partial K$ by a convex function $f_m : \mathbb{R}^{n-1}
\rightarrow \mathbb{R}$, such that $(t_1, \dots, t_{n-1})
\rightarrow \big(t_1, \dots, t_{n-1}, f_m(t_1, \dots,
t_{n-1})\big)\in \partial K$. \vskip 2mm As $\partial K$ is $C^3$,
by Taylor's theorem there exists $s \in P(K)$ such that for all
$t=(t_1, \dots, t_{n-1}) \in P(K)$, the orthogonal projection of
$K$ onto $\mathbb{R}^{n-1}$,
\begin{equation*}\label{Taylor}
f_m(t)=\frac{1}{2} \langle t, A t \rangle +\sum_{ |\alpha|=3}
\frac{D^\alpha f_m(s)}{\alpha !} t^\alpha,
\end{equation*}
where $t^\alpha=t_1^{\alpha_1} \dots  t_{n-1}^{\alpha_{n-1}}$ and
$A= \left(\frac{\partial ^2 f_m}{\partial t_i \partial
t_j}(0)\right)_{i,j = 1}^{n-1}$ is the Hessian of $f_m$ at $0$.
Clearly,  $A=\frac{1}{r_m}Id_{n-1}$ with $Id_{n-1}$ the identity
matrix.  $K\in C^3$ also implies that $|D^\alpha f_m(t)|\leq D$
for some $D>0$, for all $\a$ with $|\a|=3$, and all $t\in P(K)$.
Therefore,
\begin{equation}\label{Taylor:1}
f_m(t)\leq \frac{1}{2} \langle t, A t \rangle + \frac{(n-1)^3}{6} D
\|t\|^3=\frac{1}{2}\frac{\|t\|^2}{r_m}+ \frac{(n-1)^3}{6}
D\|t\|^3.
\end{equation}
\par
\noindent
$\partial B^n_2\big(x_m-\bar{r}_m,
\bar{r}_m\big)$ is described in our chosen coordinate system by a function $g_m: \mathbb{R}^{n-1}
\rightarrow \mathbb{R}$ - for $t_n=g_m(t) \leq \bar{r}_m $ -
$$t_n=g_m(t) =\bar{r}_m
\left[1-\sqrt{1-\left(\frac{\|t\|}{\bar{r}_m}\right)^2}\right], \
\ \mbox{if }\ \ t_n\leq \bar{r}_m.
$$ As $\frac{b}{2}\leq 1-\sqrt{1-b}$ for all $b\in [0,1]$, we get for all $t$ with $\|t\|\leq \bar{r}_m$
\begin{equation}\label{Taylor:2}
t_n=g_m(t) \geq \frac{1}{2} \frac{\|t\|^2}{(1-a)r_m}.
\end{equation}
\par
 \noindent
 Thus by  (\ref{Taylor:1}) and (\ref{Taylor:2}),  for \begin{equation*}\|t\|\leq
t_{a,m}= \min\left\{\bar{r}_m, \frac{3a}{ D \bar{r}_m
(n-1)^3}\right\},\end{equation*} one has $t_n= g_m(t) \geq
f_m(t)$. We then let
$$\Delta _{a, m}=\bar{r}_m -\sqrt{\bar{r}_m^2-t_{a,m}^2},$$ and hence
for all $\Delta \leq \Delta _{a,m}$, condition (\ref{approx:1})
holds true. We further let \begin{equation}\label{delta:m:1}
\d_{m,1}=|B^{n-1}_2|\int _{\bar{r}_m-\Delta_{a,m}}^{\bar{r}_m}
\big(\bar{ r}_m^2-y^2\big)^{\frac{n-1}{2}}\,dy.
\end{equation}
\par
Denote by $\Delta_{\bar{r}_m}$  the height of a cap of
$\bar{r}_m B^n_2$ of volume exactly $\d$. Recall
$x_{m,\d}=\partial K_\d \cap[0,x_m]$, and $$ \Delta_{x_m} =
\left\langle \frac{x_m}{\|x_m\|}, N_K(x_m)\right \rangle  \  \|x_m
- x_{m,\d}\|.
$$
By (\ref{approx:1}) and (\ref{delta:m:1}), one has, for all
$\d\leq \d_{m,1}$,
\begin{equation*}
\Delta_{x_m} \leq \Delta_{\bar{r}_m}.
\end{equation*}
Moreover, as $x$ and $x_{m,\delta}$ are colinear,
\begin{equation*}
1-\left( \frac{ \|x_ {m, \d}  \|}{\|x_m\|}\right)^n \leq \ n\
\frac{\|x_m-x_{m, \d}\|}{\|x_m\|}.
\end{equation*}
\par
\noindent
Thus
\begin{eqnarray}
\frac{c_n }{ \d^{\frac{2}{n+1}}}\ \frac{\|x_m-x_{m,
\d}\|}{\|x_m\|} = \frac{c_n\ \D _{x_m}}{ \d^{\frac{2}{n+1}} \
\langle x_m, N_K(x_m) \rangle} \label{above0} \leq \frac{c_n\ \D
_{\bar{r}_m}}{ \d^{\frac{2}{n+1}} \ \langle x_m, N_K(x_m)
\rangle}. \label{above}
\end{eqnarray}

\par
\noindent
Let  $\d < \frac{|B^n_2|\bar{r}_m^n}{2}$. Then
$\D_{\bar{r}_m}<\bar{r}_m$ and by definition of
$\Delta_{\bar{r}_m}$,
\begin{eqnarray*}
 \d &=& |B^{n-1}_2| \ \int _{\bar{r}_m-\D _{\bar{r}_m}} ^{\bar{r}_m}
(\bar{r}_m-y)^{\frac{n-1}{2}}\ (\bar{r}_m+y)^{\frac{n-1}{2}}\,dy \\
&\geq& |B^{n-1}_2| \ (2\bar{r}_m-\D_{\bar{r}_m})^{\frac{n-1}{2}}
\int _{\bar{r}_m-\D _{\bar{r}_m}} ^{\bar{r}_m}
(\bar{r}_m-y)^{\frac{n-1}{2}}\,dy\\&=&2^\frac{n+1}{2}  \
\bar{r}_m^\frac{n-1}{2} \  \frac{|B^{n-1}_2|}{n+1} \
\left(1-\frac{\D_{\bar{r}_m}}{2\bar{r}_m}\right)^{\frac{n-1}{2}}
\D_{\bar{r}_m}^{\frac{n+1}{2}}.
\end{eqnarray*}
In particular, using $\left(1-\frac{\D_{a,m}}{2\bar{r}_m} \right)
>\frac{1}{2}$ and $\Delta _{a, m}=\bar{r}_m \left(1
-\sqrt{1-\frac{t_{a,m}^2}{\bar{r}_m^2}}\right) \geq
\frac{t_{a,m}^2}{2 \bar{r}_m}$,
\begin{eqnarray}\label{delta:m:0}
 \d_{m,1} \geq
 \frac{t_{a,m}^{n+1} \   |B^{n-1}_2|}{2^\frac{n-1}{2} \  (n+1) \  \bar{r}_m} = \delta_m.
 \end{eqnarray}
And
\begin{eqnarray*}
 \d &=& |B^{n-1}_2| \ \int _{\bar{r}_m-\D _{\bar{r}_m}} ^{\bar{r}_m}
(\bar{r}_m-y)^{\frac{n-1}{2}}\ (\bar{r}_m+y)^{\frac{n-1}{2}}\,dy\\
&\leq& |B^{n-1}_2| \ 2^{\frac{n-1}{2}}\  \bar{r}_m^\frac{n-1}{2}
\int _{\bar{r}_m-\D _{\bar{r}_m}} ^{\bar{r}_m}
(\bar{r}_m-y)^{\frac{n-1}{2}}\,dy =2^\frac{n+1}{2}  \
\bar{r}_m^\frac{n-1}{2} \  \frac{|B^{n-1}_2|}{n+1} \
 \D_{\bar{r}_m}^{\frac{n+1}{2}}.
\end{eqnarray*}
Hence for all $\d < \frac{\bar{r}_m^n |B^n_2|}{2}$,
\begin{equation} \label{r0}
\frac {\bar{r}_m^{-\frac{n-1}{n+1}}}{c_n}\   \leq  \
\frac{\D_{\bar{r}_m}} {\d^{\frac{2}{n+1}}} \ \leq \ \frac
{\bar{r}_m^{-\frac{n-1}{n+1}}}{c_n}
\left(1-\frac{\D_{\bar{r}_m}}{2\bar{r}_m}\right)^{-\frac{n-1}{n+1}}.
\end{equation}

\noindent From  (\ref{above}), (\ref{delta:m:0}),  and (\ref{r0})
we get for all $\d < \min \{\frac{|B^n_2|\bar{r}_m^n}{2}, \d
_{m}\}$,
\begin{eqnarray*}
\frac{c_n }{n\ \d^{\frac{2}{n+1}}}\left[1-\left(\frac{\|{x}_{m,
\d}\|}{\|{x}_m\|}\right)^n\right] &\leq& \frac {\bar{r}_m^{-\frac{n-1}{n+1}}}{\langle {x}_m, N_K({x}_m)\rangle}
\left(1-\frac{\D_{\bar{r}_m}}{2\bar{r}_m}\right)^{-\frac{n-1}{n+1}} \\
&= & \left(1-a\right)^{-\frac{n-1}{n+1}}\frac {\kappa_K({x}_m)^\frac{1}{n+1}}{\langle {x}_m, N_K({x}_m)\rangle}
\left(1-\frac{\D_{\bar{r}_m}}{2\bar{r}_m}\right)^{-\frac{n-1}{n+1}} \\
&= & {T}_m \ \left(1-a\right)^{-\frac{n-1}{n+1}}\
\left(1-\frac{\D_{\bar{r}_m}}{2\bar{r}_m}\right)^{-\frac{n-1}{n+1}}.
\end{eqnarray*}
\par
\noindent Hence for all $\d \leq \min\{\d_2,
\frac{|B^n_2|\bar{r}_m^n}{2}, \d _{m}\}$, (\ref{xm:estimate}) of
the previous section will hold, if we choose $\d_2$ so that
\begin{equation}\label{d2/1}
\left(1-\frac{\D_{\bar{r}_m}}{2\bar{r}_m}\right)^{-\frac{n-1}{n+1}}
\leq \left(1-a\right)^{\frac{n-1}{n+1}} \ \left(\frac{\t
+1}{2}\right).
\end{equation}
By (\ref{r0}) and (\ref{d2/1}),
$$
\D_{\bar{r}_m} \leq  \d^{\frac{2}{n+1}} \  \frac
{\bar{r}_m^{-\frac{n-1}{n+1}}}{c_n}
\left(1-\frac{\D_{\bar{r}_m}}{2\bar{r}_m}\right)^{-\frac{n-1}{n+1}}
\leq \d^{\frac{2}{n+1}}  \  \frac
{\bar{r}_m^{-\frac{n-1}{n+1}}}{c_n} \
\left(1-a\right)^{\frac{n-1}{n+1}} \ \left(\frac{\t +1}{2}\right).
$$
Thus, to have (\ref{d2/1}),  it is enough to have that
$$
\left(1- \frac{\d^\frac{2}{n+1}}{c_n \bar{r}_m^\frac{2n}{n+1}}\
\left(\frac{{\tau}+1}{4}\right) \left(1-a\right)^{\frac{n-1}{n+1}}
\right)^{-1} \leq  \left(1-a\right) \left(\frac{{\tau} +1}{2}
\right)^\frac{n+1}{n-1}.
$$
Hence, we can let
\begin{equation}\label{d2/1final}
\d_2= 2^{\frac{3(n+1)}{2}}\  \left(\frac{1-a}{1+\tau}\right)^\frac{n+1}{2}\ {r}_m^n \  \frac{|B^{n-1}_2|}{n+1} \ \left[
1-(1-a)^{-1}\left(\frac{2}{{\tau}+1}\right)^\frac{n+1}{n-1}
\right]^\frac{n+1}{2}.
\end{equation}

\vskip 2mm Now we consider $x_M$. We let
$r_M=\k_K(x_M)^{\frac{1}{1-n}}$, $\bar{r}_M=(1-a)r_M$ and
$\bar{R}_M=(1+a)r_M.$ Assume now that $x_M=0$, that
$N_K(x_M)=-e_n$ and  that the other $(n-1)$ axes of the
approximating ellipsoid coincide with the remaining $(n-1)$
coordinate axes. Locally  we can then describe $\partial K$ by a
convex function $f_M: \mathbb{R}^{n-1} \rightarrow \mathbb{R}$,
such that $(t_1, \dots, t_{n-1}) \rightarrow \big(t_1, \dots,
t_{n-1}, f_M(t_1, \dots, t_{n-1})\big)\in \partial K$.

\par
\noindent
For
$$t_{M,1}=\min\left\{\bar{r}_M, \frac{3a}{ D
\bar{r}_M (n-1)^3}\right\},
$$
let \begin{equation*} \Delta_{M,1}=
\bar{r}_M-\sqrt{\bar{r}^2_M-t_{M,1}^2}.
\end{equation*}
We repeat  the previous argument  and get  for all
$\Delta \leq \Delta_{M,1}$, that
\begin{eqnarray}\label{approx:M:1}
B^n_2\left(x_M-\bar{r}_M N_K(x_M), \bar{r}_M \right)  \cap H^-\left(x_M- \Delta N_K(x_M), N_K(x_M) \right)  \nonumber \\
 \subseteq K \cap H^-\left(x_M- \Delta N_K(x_M), N_K(x_M) \right).
\end{eqnarray}

\par
\noindent Let $\xi=1+\frac{a}{2}$. Then for all $b\in [0,
\frac{4(\xi-1)}{\xi^2}]$,
$$1-\sqrt{1-b}\leq \frac{\xi}{2}b.
$$
For $g_M(t) \leq \bar{R}_M$, let
$$t_n=g_M(t) =\bar{R}_M
\left[1-\sqrt{1-\left(\frac{\|t\|}{\bar{R}_M}\right)^2}\right],
$$
be the function that describes $\partial B^n_2\left(x_M-\bar{R}_M
N_K(x_M), \bar{R}_M \right)$. For $\|(t_1 \dots,t_{n-1}) \| \leq
\frac{2(\xi-1)^\frac{1}{2}}{\xi} \ \bar{R}_M$,  we get that
\begin{equation*}
t_n= g_M(t)\leq \frac{\xi}{2} \frac{\|t\|^2}{\bar{R}_M}.
\end{equation*}
As in inequality (\ref{Taylor:1}),
\begin{equation*} f_M(t)\geq
\frac{1}{2}\frac{\|t\|^2}{r_M}-\frac{(n-1)^3}{6} D\|t\|^3.
\end{equation*}
Therefore,  for $\|t\|\leq t_{M,2}$ where
$$t_{M,2}=\min\left\{\frac{2(\xi-1)^\frac{1}{2}}{\xi} \ \bar{R}_M,
\frac{3a}{2D\bar{R}_M(n-1)^3}\right\},$$ we get that $f_M(t)\geq g_M(t)= t_n$
which implies that
\begin{eqnarray}\label{approx:M:2}
B^n_2\left(x_M-\bar{R}_M N_K(x_M), \bar{R}_M \right)  \cap H^-\left(x_M- \Delta N_K(x_M), N_K(x_M) \right)  \nonumber \\
 \supseteq K \cap H^-\left(x_M- \Delta N_K(x_M), N_K(x_M) \right),
\end{eqnarray}
for all $\Delta \leq \Delta_{M,2}$ with
$\Delta_{M,2}=\bar{R}_M-\sqrt{\bar{R}_M^2-t_{M,2}^2}.$

\vskip 2mm We let $\Delta_{a,M}=\min\{\Delta_{M,1},
\Delta_{M,2}\}$, and
\begin{equation}\label{Delta:M:M1}\d_{M,1}=|B^{n-1}_2|\int_{\bar{r}_M-\Delta_{a,
M}}^{\bar{r}_M}(\bar{r}_M^2-y^2)^{\frac{n-1}{2}}\,dy.
\end{equation}
In particular, as above,
\begin{equation}\label{delta:M:0}
\d_{M,1} \geq 2 \frac{|B^{n-1}_2|}{n+1}\ \bar{r}_M^{\frac{n-1}{2}}
\Delta_{a,M} ^\frac{n+1}{2}= \delta_M,
\end{equation} where $ \Delta_{a,M}$ can be taken as $$ \min
\left\{\frac{t_{M,1}^2}{2(1-a)r_m},
\frac{t_{M,2}^2}{2(1-a)R_M}\right\}.$$

\par \noindent
 Let $\D _{\bar{R}_M}$ be the height of a cap of
$\bar{R}_MB^n_2$ of volume exactly $\delta$. By
(\ref{approx:M:1}), (\ref{approx:M:2}),  (\ref{Delta:M:M1}), and
(\ref{delta:M:0}), for all $\delta \leq \delta _{M}$, $$\Delta
_{\bar{R}_M} \leq \D _{x_M}=\left\langle \frac{x_M}{\|x_M\|},
N_K(x_M)\right \rangle \ \|x_M - x_{M,\d}\|.$$  Recall
$x_{M,\d}=[0,x_M]\cap
\partial K_{\d}$. Similar to (\ref{r0}), for all $\d <\min\{
\frac{\bar{r}_M^n |B^n_2|}{2}, \d _{M}\}$,
\begin{equation*} \frac {\bar{R}_M^{-\frac{n-1}{n+1}}}{c_n}\   \leq  \
\frac{\D_{\bar{R}_M}} {\d^{\frac{2}{n+1}}} \ \leq \ \frac
{\bar{R}_M^{-\frac{n-1}{n+1}}}{c_n}
\left(1-\frac{\D_{\bar{R}_M}}{2\bar{R}_M}\right)^{-\frac{n-1}{n+1}}.
\end{equation*} This implies that, for
all $\d <\min\{ \frac{\bar{r}_M^n |B^n_2|}{2}, \d _{M}\}$,
\begin{eqnarray}
\frac{c_n }{n\ \d^{\frac{2}{n+1}}}\left[1-\left(\frac{\|x_{M,
\d}\|}{\|x_M\|}\right)^n\right]\!\! &=&\!\! \frac{c_n }{n\
\d^{\frac{2}{n+1}}}\left[1-\left(1-\frac{\Delta_{x_M}
}{\langle x_M, N_K(x_M)\rangle}\right)^n\right]  \nonumber\\
\!\!&\geq&\!\!\frac{c_n }{n\ \d^{\frac{2}{n+1}}}\left[1-\left(1-\frac{\Delta_{\bar{R}_M}}{\langle x_M, N_K(x_M)\rangle}\right)^n\right] \nonumber \\
\!\!&\geq&\!\! \frac{c_n }{ \d^{\frac{2}{n+1}}} \left( \!\! \frac{
\D _{\bar{R}_M}}{   \langle x_M, N_K(x_M) \rangle} - \frac{(n-1)\
\D ^2_{\bar{R}_M}}{2\  \langle x_M, N_K(x_M) \rangle^2}\!\!
\right)\nonumber \\ &\geq & \frac
{\bar{R}_M^{-\frac{n-1}{n+1}}}{\langle {x}_M, N_K({x}_M)\rangle}
\left( 1- \frac{(n-1) \D_{\bar{R}_M}}{2 \langle {x}_M,
N_K({x}_M)\rangle} \right)\nonumber\\ \hskip 37mm &\geq &\!\!
(1+a)^{-\frac{n-1}{n+1}}\!{T}_M\! \left(\! 1\!-\!\frac{(n-1)
\d^{\frac{2}{n+1}} {T}_M}{2c_n (1+a)^{\frac{n-1}{n+1}}}
\bigg[1-\frac{\D_{\bar{R}_M}}{2\bar{R}_M}\bigg]^{-\frac{n-1}{n+1}}\right)\nonumber\\
&\geq& (1+a)^{-\frac{n-1}{n+1}} {T}_M \   \left(1- \frac{(n-1) \
{T}_M\ \d^\frac{2}{n+1} }
{2^\frac{2}{n+1}c_n(1+a)^{\frac{n-1}{n+1}}}\right). \nonumber
\end{eqnarray}
The last inequality holds as  $\frac{\D_{\bar{R}_M}}{2\bar{R}_M} \leq
\frac{1}{2}.$ To have (\ref{xM:estimate}) of the previous section, it is enough to have
$$ (1+a)^{-\frac{n-1}{n+1}}   \left(1- \frac{(n-1) \
{T}_M\ \d^\frac{2}{n+1} }
{2^\frac{2}{n+1}c_n(1+a)^{\frac{n-1}{n+1}}}\right) \geq \frac{2\t
+1 }{3\t},$$ or, equivalently,  \begin{equation} \label{delta:2:2}
\d \leq
\frac{[1-(1+a)^{\frac{n-1}{n+1}}\left(\frac{2\t+1}{3\t}\right)]^{\frac{n+1}{2}}\
2^{\frac{n+3}{2}} \ |B^{n-1}_2| (1+a)^{\frac{n-1}{2}}
}{(n-1)^{\frac{n+1}{2}}T_M^{\frac{n+1}{2}}(n+1)}(=:\d
_1).\end{equation}

\vskip 2mm \noindent Now we let the threshold $\d(K)$ be
\begin{eqnarray*} \d(K)&=&\min\bigg\{\delta_0,  \delta_1, \delta_2, \d _{m},  \d _{M}, \
\frac{(1-a)^{n}\ r_m^n\ |B^n_2| }{2},  \frac{(1-a)^{n}\ r_M^n \
|B^n_2|}{2}\bigg\},
\end{eqnarray*}
where $\d_0, \d_1, \d_2, \d_m, \d_M$ are as in (\ref{delta:0:0}),
(\ref{delta:2:2}), (\ref{d2/1final}), (\ref{delta:m:0}) and
(\ref{delta:M:0}) respectively.

\vskip 2mm \noindent Elisabeth Werner\\
{\small Department of Mathematics \ \ \ \ \ \ \ \ \ \ \ \ \ \ \ \ \ \ \ Universit\'{e} de Lille 1}\\
{\small Case Western Reserve University \ \ \ \ \ \ \ \ \ \ \ \ \ UFR de Math\'{e}matique }\\
{\small Cleveland, Ohio 44106, U. S. A. \ \ \ \ \ \ \ \ \ \ \ \ \ \ \ 59655 Villeneuve d'Ascq, France}\\
{\small \tt elisabeth.werner@case.edu}\\ \\
\and Deping Ye\\
{\small Department of Mathematics}\\
{\small 202 Mathematical Sciences Bldg}\\
{\small University of Missouri}\\
{\small Columbia, MO 65211 USA}\\
{\small \tt yed@missouri.edu}

\end{document}